\documentclass[11pt]{article}

\usepackage{fancyhdr,mathptmx}
\usepackage{pictexwd}
\usepackage{amsmath}
\usepackage{amscd}
\usepackage{amssymb}
\usepackage{amsthm}
\usepackage{xspace}
\usepackage[mathscr]{eucal}
\usepackage[all,tips]{xy}
\usepackage[dvips]{graphicx}
\usepackage{fancyhdr}
\usepackage{verbatim}
\usepackage{syntonly}

\fancypagestyle{plain}

\theoremstyle{plain}

\newtheorem{thm}{Theorem}[section]
\newtheorem{cor}[thm]{Corollary}

\newenvironment{pf}
{\begin{proof}} {\end{proof}}

\DeclareMathOperator{\id}{id}

\DeclareMathOperator*{\diag}{diag}

\DeclareMathOperator{\Inv}{Inv}

\DeclareMathOperator*{\ch}{char} 

\DeclareMathOperator{\Gal}{Gal}

\DeclareMathOperator{\SL}{SL}

\DeclareMathOperator{\GL}{GL}

\DeclareMathOperator{\Ort}{O}

\DeclareMathOperator{\SO}{SO}

\DeclareMathOperator{\Uni}{U}

\DeclareMathOperator{\SU}{SU}

\DeclareMathOperator{\Sp}{Sp}

\DeclareMathOperator{\Mat}{M}

\newcommand{\al}{\alpha}

\newcommand{\be}{\beta}

\newcommand{\ga}{\gamma}
\newcommand{\Ga}{\Gamma}

\newcommand{\te}{\theta}

\newcommand{\si}{\sigma}

\newcommand{\la}{\lambda}
\newcommand{\La}{\Lambda}

\newcommand{\ol}{\overline}

\newcommand{\fn}{\footnote}

\newcommand{\wh}{\widehat}

\newcommand{\co}{\ensuremath{\colon}}

\newcommand{\innp}[1]{\left< #1 \right>}
\newcommand{\abs}[1]{\left\vert#1\right\vert}
\newcommand{\set}[1]{\left\{#1\right\}}

\newcommand{\su}{\subset}

\newcommand{\op}{\oplus}

\newcommand{\bop}{\bigoplus}

\newcommand{\smin}{\setminus}

\newcommand{\lra}{\longrightarrow}

\newcommand{\B}[1]{\ensuremath{\mathbf{#1}}}

\newcommand{\Cal}[1]{\ensuremath{\mathcal{#1}}}

\newcommand{\Hy}{\ensuremath{\B{H}}}
\newcommand{\N}{\ensuremath{\B{N}}}
\newcommand{\Q}{\ensuremath{\B{Q}}}
\newcommand{\R}{\ensuremath{\B{R}}}
\newcommand{\Z}{\ensuremath{\B{Z}}}
\newcommand{\C}{\ensuremath{\B{C}}}

\newcommand{\refS}[1]{Section~\ref{S:#1}}

\newcommand{\refT}[1]{Theorem~\ref{T:#1}}

\begin{document}

%---------------------------------------------------------------------------
%---------------------------------------------------------------------------

\title{\textbf{Finite subgroups of \\ arithmetic lattices in $\Uni(2,1)$}}
\author{D. B. McReynolds\fn{Supported by a Continuing Education fellowship}}
\maketitle
%---------------------------------------------------------------------------
%---------------------------------------------------------------------------

\begin{abstract}
\noindent The principle result of this article is the determination of the possible finite subgroups of arithmetic lattices in $\Uni(2,1)$.
\end{abstract} 

%---------------------------------------------------------------------------
%---------------------------------------------------------------------------

%---------------------------------------------------------------------------
%---------------------------------------------------------------------------
\section{Introduction}

\noindent The study of finite subgroups in linear algebraic groups and their discrete subgroups is relatively old. One of the oldest results is due to Jordan \cite{Jordan1878} who gave an upper bound on the order of a finite subgroup of $\GL(n;\Z)$ (see also Boothby--Wang \cite{BoothbyWang65}, Brauer--Feit \cite{BrauerFeit66}, Feit \cite{Feit97}, and Weisfeiler \cite{Weisfeiler84A}). Any such subgroup is conjugate into the orthogonal group $\Ort(n)$ and when generated by reflections, admit a special presentation by work of Coxeter---see \cite{BensonGrove70}. Friedland \cite{Friedland97} generalized Jordan's theorem by giving an upper bound for the order of a finite subgroup of $\GL(n;\Q)$---see also \cite{AljadeffSonn98}. Such a subgroup is visibly a subgroup of the central simple $\Q$--algebra $\Mat(n;\Q)$. In this vein, Vign\'{e}ras \cite{Vigneras80} classified finite subgroups of $\Q$--defined quaternion algebras; below we discuss Amitsur's work on finite subgroups of division algebras with characteristic zero central fields. For arithmetic lattices in other Lie groups, Ratcliffe and Tschantz \cite{RatcliffeTschantz99} investigated finite subgroups of $\Ort(n,1;\Z)$, obtaining upper bounds on the order of a possible finite subgroup and sizes of $p$--subgroups.\smallskip\smallskip

\noindent More germane to this article are the finite subgroups of unitary groups. In \cite{ShephardTodd54} (see also \cite[p. 98]{Coxeter74} or \cite[p. 57]{DuVal64}), the finite subgroups of $\Uni(2)$ and the finite complex reflection groups were classified. For those reflection groups in $\Uni(2) \times \Uni(1)$, Falbel--Paupert \cite{FalbelPaupert04} constructed fundamental domains for the induced action of these groups on $\Hy_\C^2$. In addition, complex reflections groups have also been thoroughly investigated by Brou{\'e}--Malle--Rouquier \cite{BMR98} who, among other things, computed invariants for these groups.\smallskip\smallskip

\noindent The arithmetic lattices in $\Uni(2,1)$ come in two flavors. The better known are those lattices commensurable with $\Uni(H;\Cal{O}_E)$, where $H$ is an admissible hermitian form over a CM field $E/F$ (see \S 2 for terminology). The other family of lattices are derived from cyclic central division algebras over $E/F$ equipped with an involution of second kind (see \cite[p. 83]{PlatonovRapinchuk94}). These lattices are called \emph{first} and \emph{second type}, respectively.\smallskip\smallskip

\noindent The purpose of this short article is to classify the possible finite subgroups of arithmetic lattices of first and second type in $\Uni(2,1)$. Any such subgroup is conjugate in $\GL(3;\C)$ into $\Uni(2) \times \Uni(1)$. In particular, this is a more refined form of the determination of the finite subgroups of $\Uni(2) \times \Uni(1)$ and can be viewed as a precise form of Jordan's theorem for these classes of lattices. Our first result is the realization of any finite subgroup of $\Uni(2) \times \Uni(1)$ as a subgroup of an arithmetic lattice of first type.\smallskip\smallskip

\begin{thm}\label{T:FiniteInForms} If $G$ is a finite subgroup of $\Uni(2) \times \Uni(1)$, then there exists an arithmetic lattice $\Ga$ of first type in $\Uni(2,1)$ which contains a subgroup isomorphic to $G$.
\end{thm}

\noindent For arithmetic lattices of the second type, the story is quite different.\smallskip\smallskip

\begin{thm}\label{T:CST} If $\Ga$ is a arithmetic lattice of second type in $\Uni(2,1)$ with finite subgroup $G$, then $G$ is cyclic.
\end{thm}

\noindent The proof of \refT{CST} follows from a more general result for finite subgroups of arithmetic lattices of second type in $\Uni(p-1,1)$ where $p$ is an odd prime. Using local theory, for a fixed commensurability class of lattices $[\Ga]$ of second type, the possible cyclic subgroups can be determined---see \cite[Ch. 12.5]{MaclachlanReid03} for a treatment in the Fuchsian and Kleinian setting.\smallskip\smallskip

\noindent That these two classes of lattices are substantially different is well known. By work of Kazhdan \cite{Kazhdan77}, lattices of first type virtually surject $\Z$, while Rogawski \cite{Rogawski90} showed congruence subgroups of principal arithmetic lattices of second type in $\Uni(2,1)$ do not. In addition, Reznikov \cite{Reznikov95} (see also Klingler \cite{Klingler03}) proved a superrigidity-type theorem for representations of congruence subgroups of principal arithmetic lattices of second type into $\GL(3)$. Recently, Stover \cite{Stover05A} proved congruence subgroups of principal arithmetic lattices of second type cannot split as nontrivial amalgamated products. This class includes all the known so-called fake projective planes---see for instance Mumford \cite{Mumford79}.\smallskip\smallskip

\noindent In the final section of this article, we briefly discuss this general problem for lattices in Lie groups. In the case of arithmetic lattices in the real classical groups, much of this article extends without fuss.

%---------------------------------------------------------------------------
\subsection*{Acknowledgements}

\noindent I would like to thank my advisor Alan Reid for all his help and Elisha Falbel for many invaluable comments and suggestions. A debt is also owed to Matthew Stover for carefully reading early drafts of this article.

%---------------------------------------------------------------------------
%---------------------------------------------------------------------------
\section{Preliminaries}

\noindent We briefly recall requisite background material here and establish some notation to be used in the remainder of this article. \smallskip\smallskip

\noindent By a \emph{totally real number field} (resp. \emph{totally imaginary}), we mean an finite field extension $F$ of $\Q$ such that each field embedding of $F$ into $\C$ is real (resp. complex). By a \emph{CM field}, we mean a quadratic extension $E$ of $F$ for which $E$ is totally imaginary and $F$ is totally real. We list the distinct embeddings of $E \lra \C$ by $\tau_1,\dots,\tau_q$ which are taken up to field automorphisms of $\C$. For each embedding $\tau_j$ of $E$, we have an associated real embedding $\si_j$ of $F$ given by restriction and every real embedding of $F$ uniquely arises this way up to field automorphisms of $\R$. We fix once and for all an embedding $\tau_1\co E \lra \C$ and hence consider $E \su \C$ and $F \su \R$ via $\tau_1$ and $\si_1$. \smallskip\smallskip

\noindent For a Galois extension $L/E$ with Galois group $\Gal(L/E)$, the \emph{norm} of an element $\ga \in L$ is defined to be
\[ N_{L/E}(\al) = \prod_{g \in \Gal(L/E)} g(\al). \]
It follows at once that $N_{L/E}(\al) \in E$. Less obvious is that for a cyclic extension $L/E$ of degree $n$, the index of $N_{L/E}(L^\times)$ in $E^\times$ is $n$ with quotient $E^\times/N_{L/E}(L^\times)$ isomorphic to $\Z/n\Z$; this is a deep result from class field theory.\smallskip\smallskip

\noindent In the sequel, we will make repeated use of the following well known theorem from algebraic number theory.

\begin{thm}[Weak approximation theorem]\label{T:WeakApprox2}
Let $K$ be a totally real number field with real embeddings $\si_1,\dots,\si_s$. For any $p,q \in \N$ with $p+q=s$, there exists $\la \in K$ and embeddings $\si_{j_1},\dots,\si_{j_p}$ such that $\si_k(\al)>0$ if and only if $k=j_\ell$ for $\ell=1,\dots,p$.
\end{thm}

%---------------------------------------------------------------------------
%---------------------------------------------------------------------------
\section{Arithmetic lattices in $\Uni(n,1)$}

%---------------------------------------------------------------------------
\subsection{Arithmetic lattices of first type}

\noindent By an \emph{$E$--defined hermitian matrix} $H$ we mean a hermitian matrix in $\GL(n;E)$. Via Gram-Schmidt, the matrix $H$ can be diagonalized over $\C$ and possesses real nonzero eigenvalues. We denote the number of positive (resp. negative) eigenvalues of $H$ by $e_+(H)$ (resp. $e_-(H)$). We call the pair $(e_+(H),e_-(H))$ the \emph{signature pair} for $H$ and the number $\si(H) = \abs{e_+(H) - e_-(H)}$ the \emph{signature}. For any matrix $X \in \GL(n;E)$, each embedding $\tau_j$ of $E$ yields a new matrix $^{\tau_j}X$ by applying $\tau_j$ to the coefficients of $X$. We say the pair $(E/F,H)$ is \emph{admissible} if
\[ \si(~^{\tau_j}H) = \begin{cases} n-2, & \tau_j = \tau_1, \\ n, & \text{otherwise}. \end{cases} \]
For an admissible pair $(E/F,H)$, the group $\Uni(H;\Cal{O}_E)$ is a lattice in $\Uni(n-1,1)$ for any selection of a Lie isomorphism between $\Uni(H)$ and $\Uni(n-1,1)$. We call a lattice constructed this way a \emph{principal arithmetic lattice of first type}. Any subgroup $\La$ of $\Uni(n-1,1)$ commensurable with a principal arithmetic lattice of first type $\Ga$, i.e. $\La \cap \Ga$ is finite index in both $\La$ and $\Ga$, is called an \emph{arithmetic lattice of first type}. By the weak approximation theorem, admissible pairs $(E/F,H)$ exist for any CM field $E/F$. \smallskip\smallskip

%---------------------------------------------------------------------------
\subsection{Equivalence of hermitian forms}

\noindent For any hermitian matrix $H$ defined on an $E$--defined vector space $V$, we associate to $H$ a triple invariant. First is the $E$--dimension of $V$ which we denote by $\dim V$. Second, for each embedding $\tau_j$ of $E$, $^{\tau_j}H$ has an associated signature $\si_{\tau_j}(H)$, and we denote by $\si(H)$ the set $\set{\si_{\tau_j}(H)}_j$. Finally, by selecting an $E$--basis $\Cal{B}$ for $V$, we associate to $H$ the determinant $\det_{\Cal{B}} H$ of the associated matrix for $H$ in the basis $\Cal{B}$. This is not well defined as an element of $F^\times$, as changing the basis $\Cal{B}$ can change $\det_{\Cal{B}}(H)$. However, as an element of $F^\times/N_{E/F}(E^\times)$, it is well defined and denoted by $\det H$. We denote the triple $(\dim V, \si(H),\det H)$ by $\Inv(H)$ and call it the associated invariant for $H$. \smallskip\smallskip

\noindent We say a pair of $E$--defined hermitian matrices $H_1$ and $H_2$ on $V$ are \emph{equivalent} if the unitary groups $\Uni(H_1)$ and $\Uni(H_2)$ are conjugate in $\GL(V;E)$---alternatively, the associated hermitian forms are similar. In this case, the groups $\Uni(H_1;\Cal{O}_E)$ and $\Uni(H_2;\Cal{O}_E)$ are commensurable in the wide sense. Namely, a $\GL(V)$--conjugate of $\Uni(H_1;\Cal{O}_E)$ is commensurable with $\Uni(H_2;\Cal{O}_E)$. The following result can be found in \cite[Cor. 6.6, p. 376]{Scharlau85}.\smallskip\smallskip

\begin{thm}\label{T:FormClassification} Two $E$--defined hermitian matrices $H_1$ and $H_2$ on $V$ are equivalent if and only if $\Inv(H_1)=\Inv(H_2)$.
\end{thm}

%---------------------------------------------------------------------------
\subsection{Arithmetic lattices of second type}

\noindent This article does not require an in depth discussion on the construction of arithmetic lattices of second type. For our purposes we need only know any lattice $\Ga$ of second type possesses a faithful representation into $A^\times$ for some cyclic division algebra $A$ whose central field is characteristic zero. When $\Ga<\Uni(n-1,1)$, the degree of the algebra $A$ over its central field is $n$. Nevertheless, for completeness, we briefly describe the construction of these lattices in $\Uni(2,1)$.\smallskip\smallskip 

\noindent Let $E/F$ be a CM field with Galois involution $\te$ and $L/E$ a cyclic Galois extension of degree three with nontrivial Galois automorphism $\tau$. The field $L$ possesses a unique totally real subfield $K$ and these four fields fit into the diagram
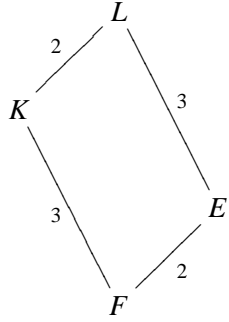
\[ \xymatrix{ & L \ar@{-}[ld]_2 \ar@{-}[rdd]^3 & \\ K \ar@{-}[rdd]_3 & & \\ & & E \ar@{-}[ld]^2 \\ & F & } \]
For $\al \in E^\times \smin N_{L/E}(L^\times)$ such that there exists $\be \in K$ with $N_{E/F}(\al) = N_{K/F}(\be)$, we define an $E$--algebra $A$ by
\[ A(L/E,\tau,\al) = \set{\be_0 + \be_1 X + \be_2X^2~:~\be_j \in L}, \]
subject to the relations $X^3 = \al$, $X\be = \tau(\be)X$, for $\be \in L$. With our selection of $\al$, the algebra $A$ is a division algebra with center $E$ and is degree three over $E$; $A$ is a 9--dimensional $E$--vector space. In addition, $A \otimes_E L = \Mat(3;L)$, $A$ admits an involution $\star$ such that $\star_{|E} = \te$, and the extension of $\star$ to $A\otimes_E L$ is complex transposition. Such an involution is called \emph{second kind}. To build a lattice in $\Uni(2,1)$, we select a $\star$--hermitian element $h \in A^\times$, i.e. an $h$ such that $h^\star = h$, and define
\[ \Uni(h;A) = \set{x \in A~:~ h^{-1}x^\star h x = \la,~\la \in E, \quad N_{A/E}(\la) = 1}. \]
By a \emph{$\Cal{O}_E$--order} $\Cal{O}$ of $A$, we mean a finitely generated $\Cal{O}_E$--submodule of $A$ such that $\Cal{O} \otimes_{\Cal{O}_E} E = A$. Selecting a $\Cal{O}_E$--order $\Cal{O}$ of $A$, the subgroup $\Uni(h;\Cal{O}) = \Uni(h;A) \cap \Cal{O}$ is a lattice in $\Uni(2,1)$ so long as the image of $h$ under the isomorphism $A \otimes_E L \lra \Mat(3;L)$ has signature pair $(2,1)$ (and is anisotropic at the other embeddings $\tau_j$). Any lattice commensurable with a lattice constructed in this way is called of \emph{second type}.\smallskip\smallskip

\noindent Perhaps the most well known second type lattice is the fundamental group of the complex hyperbolic 2--manifold constructed by Mumford in \cite{Mumford79} known by many as Mumford's fake $\B{CP}^2$. Other examples are unitary Shimura varieties---see \cite{HarrisTaylor01} or \cite{Rogawski90}. 

%---------------------------------------------------------------------------
%---------------------------------------------------------------------------
\section{Finite subgroups in first type lattices}

\noindent Having reviewed the requisite material, we begin the main body of this article. In this section, we investigate the possible finite subgroups of arithmetic lattices in $\Uni(2,1)$ and more generally $\Uni(n,1)$. The results of this section are neither difficult nor surprising. Their inclusion is primarily for comparative purposes with the analogous results for second type lattices. We are unaware of existing results of this flavor for these lattices---see \refS{Final} for more on finite subgroups of other arithmetic groups---but would not be surprised if some or all of the results of this section were known.

%---------------------------------------------------------------------------
\subsection{Proof of \refT{FiniteInForms}}

\noindent We now prove \refT{FiniteInForms}.\smallskip\smallskip

\begin{pf}[Proof of \refT{FiniteInForms}] For any finite subgroup $G$ of $\Uni(2)\times\Uni(1)$, there exists a CM field $E/F$ and a faithful representation $\rho\co G \lra \Uni(2;\Cal{O}_E) \times \Uni(1;\Cal{O}_E)$. This can be verified using the classification finite subgroups of $\Uni(2)$ given in \cite[p. 98]{Coxeter74} or \cite[p. 57]{DuVal64}. Indeed, each finite subgroup $G_0$ of $\Uni(2)$ has a faithful representation into $\Uni(2;\Z(\zeta_r))$ for some primitive $r$th root of unity. From this we get a faithful representation of $G$ into $\Uni(2;\Z(\zeta_{r^\prime})) \times \Uni(1;\Z(\zeta_{r^\prime}))$, where $\zeta_{r^\prime}$ is a primitive $r^\prime$th root of unity which is determined by the finite cyclic subgroup of $\Uni(1)$ and $\zeta_r$.
For any $\al \in F^\times$, define the hermitian matrix
\[ H_\al = \begin{pmatrix} 1 & 0 & 0 \\ 0 & 1 & 0 \\ 0 & 0 & \al \end{pmatrix}. \]
As $G$ is contained in $\Uni(2) \times \Uni(1)$, its follows $G$ is contained in $\Uni(H_\al;\Cal{O}_E)$. By the weak approximation theorem, we can select $\al \in F^\times$ such that $\al<0$ and for each nontrivial embedding $\si_\ell$ of $F$, $\si_\ell(\al)>0$. For this selection of $\al$, $(E/F,H_\al)$ is an admissible pair, and so $\Uni(H_\al;\Cal{O}_E)$ is an arithmetic lattice of first type containing a subgroup isomorphic to $G$.
\end{pf}

\noindent For the CM field $E/F$, the form $H_\al$ is unique up to equivalence. Consequently, the wide commensurability class of $\Ga$ is unique, up to changing the representation of $G$ into $\Uni(2) \times \Uni(1)$ and changing the CM field $E/F$. By taking any CM field $E^\prime/F^\prime$ with $E \su E^\prime$ and $F\su F^\prime$, the construction above produces a new commensurability class $[\La]$ which contains a representative $\La$ with a finite subgroup isomorphic to $G$. In particular, there exist infinitely many distinct commensurability classes of lattices of first type in $\Uni(2,1)$ which contain a representative having a finite subgroup isomorphic to $G$.\smallskip\smallskip

\begin{cor} For every finite subgroup $G$ of $\Uni(2) \times \Uni(1)$, there exist infinitely many distinct commensurability classes $[\Ga_j]$ of lattices of first type in $\Uni(2,1)$ which contain a representative $\Ga_j$ and a subgroup $G_j<\Ga_j$ isomorphic to $G$.
\end{cor}

%---------------------------------------------------------------------------
\subsection{Complex reflection groups}

\noindent It is not difficult to see the proof of \refT{FiniteInForms} works for any finite subgroup $G$ of $\Uni(n) \times \Uni(1)$ which is conjugate into $\GL(n;\Cal{O}_E)\times \GL(1;\Cal{O}_E)$ for some CM field $E/F$ to produce a faithful representation of $G$ into an arithmetic complex hyperbolic lattice $\La$ of first type in $\Uni(n,1)$. This likely can be used to realize many finite complex reflection groups (see the appendix in \cite{BMR98}) as finite subgroups of arithmetic complex hyperbolic lattices; here the dimension of the lattice is one less than the dimension of the reflection group. We thank the referee for bringing this to our attention.

%---------------------------------------------------------------------------
\subsection{Representing general finite subgroups in lattices of first type}

\noindent The following result serves to further illustrate the difference between lattices of first and second type. \smallskip\smallskip

\begin{thm}\label{T:AllFiniteAllClasses} For every finite group $G$, there exists a positive integer $n_G$ such that the following holds: For every $n\geq n_G$ and every commensurability class $[\Ga]$ of first type lattices in $\Uni(n,1)$, there exist a representative $\Ga$ which has a subgroup isomorphic to $G$.
\end{thm}

\begin{pf} To begin, every finite group $G$ admits a faithful representation (for some $m$ which depends on $G$)
\[ \rho\co G \lra \GL(m;\Z). \]
For instance, one can take the left regular representation of $G$. Select $\rho$ among all the representations of $G$ so that $m$ is minimal and denote this minimal $m$ by $n_G$. For any positive definite bilinear form $B$ defined on $\Z^{n_G}$, let $B_G$ be the $\rho(G)$--average defined by
\[ B_G(x,y) = \sum_{g \in G} B(gx,gy). \]
By construction, $B_G$ is a $\Z$--defined, positive definite, $\rho(G)$--invariant bilinear form. For any imaginary quadratic extension $E/\Q$, we can view $B_G$ as a hermitian form defined over $E$ and we denote the resulting $E$--defined hermitian form by $H_G$, which is clearly $\rho(G)$--invariant. With the form $H_G$, define yet another hermitian form by
\[ H_{G,-1} = H_G \op -I_1. \]
As $E$ has only one distinct complex embedding, $(E/\Q,H_{G,-1})$ is an admissible pair, and so $\Uni(H_{G,-1};\Cal{O}_E)$ is an arithmetic lattice of first type. By construction,
\[ \rho\op \id_1\co G \lra \Uni(H_{G,-1};\Cal{O}_E) \]
is a faithful representation. For every $n \geq n_G$, we extend $H_{G,-1}$ by an identity block $I_{n-n_G}$ which yields the hermitian form
\[ H_{G,-1,n} = H_G \op I_{n-n_G} \op -I_1. \]
As before, the pair $(E/\Q,H_{G,-1,n})$ is admissible and we have the faithful representation
\[ \rho\op \id_{n-n_G+1}\co G \lra \Uni(H_{G,-1,n};\Cal{O}_E). \]

\noindent By Godement's compactness theorem, every arithmetic lattice of second type in $\Uni(n,1)$ is cocompact. However, by work of Kneser, the lattices above are not as the forms $H_{G,-1}$ are isotropic when $n_G>1$. To produce cocompact lattices of first type which contain $G$, by Godement's compactness criterion, it suffices to produce lattices with associated admissible pair $(E/F,H)$ where $F \ne \Q$. For any CM field $E/F$, replace the form $H_{G,-1}$ with the form
\[ H_{G,\al} = H_G \op \al I_1 \]
where $\al < 0$ and for all $\si_\ell \ne \si_1$, $\si_\ell(\al)>0$. For this selection of $\al$, the pair $(E/F,H_{G,\al})$ is admissible and we have the faithful representation
\[ \rho\op \id_1\co G \lra \Uni(H_{G,\al};\Cal{O}_E). \]
This is extended to all $n\geq n_G$ in an identical manner. For even $n_G$, every lattice of first type is commensurable (in the wide sense) with a lattice constructed above. This is an immediate consequence of the fact the form $H_{G,\al}$ is unique up to equivalence. Otherwise, we must modify the form $H_{G,\al}$. To this end, there are two equivalence classes of admissible hermitian forms of dimension $n_G+1$ over $E/F$ and the two classes are parameterized by the determinant viewed as elements of $F^\times/N_{E/F}(E^\times)$. To obtain a $\rho(G)$--invariant admissible hermitian form $\wh{H}$ with $\det \wh{H}\ne \det H_{G,\al} \mod N_{E/F}(E^\times)$, we proceed as follows. If $H^\prime$ is a representative of the admissible equivalence class over $E/F$ such that
\[ \det H^\prime \ne \det H_{G,\al} \mod N_{E/F}(E^\times), \]
by changing the representative, we can assume $H^\prime$ is diagonal. If
\[ H^\prime = \diag(\al_1,\dots,\al_{n_G+1}), \]
set
\[ \be = \prod_{j=1}^{n_G+1} \al_j. \]
Since $H^\prime$ is admissible, there exists $\al_j$ such that $\al_j<0$ and for all nontrivial $\si_\ell$, $\si_\ell(\al_j)>0$. For the remaining $\al_k \ne \al_j$, for all $\si_\ell$, $\si_\ell(\al_k)>0$. In particular, $\be<0$ and for all nontrivial $\si_\ell$, $\si_\ell(\be)>0$. Therefore, the form
\[ \wh{H} = H_G \op \be I_1 \]
is an admissible form over $E/F$ which is $\rho(G)$--invariant. By construction, $\det(\wh{H}) = \det(H^\prime)$ and so $\wh{H}$ is a representative of the other admissible class over $E/F$. \end{pf}

\noindent Note we must modify $H_{G,\al,n}$ in an identical way when $n$ is odd in order to obtain representatives in both wide commensurability classes which contain a subgroup isomorphic to $G$.\smallskip\smallskip

%---------------------------------------------------------------------------
%---------------------------------------------------------------------------
\section{Proof of \refT{CST}}

\noindent The following result implies \refT{CST} when specialized to $p=3$.\smallskip\smallskip

\begin{thm}\label{T:PCST} If $\Ga$ is a arithmetic lattice of second type in $\Uni(p-1,1)$, where $p$ is an odd prime and $G$ is a finite subgroup of $\Ga$, then $G$ is cyclic.
\end{thm}

\noindent The remainder of this article is devoted to proof of this result. Briefly, we first reduce the possibilities for the finite group $G$ by work of Amitsur \cite{Amitsur55}. Using the representation theory for the remaining possible finite groups together with the fact $\Ga$ is contained in $\Uni(2,1)$, we deduce the group $G$ must be cyclic.\smallskip\smallskip

%---------------------------------------------------------------------------
\subsection{Amitsur's $D$--groups}

\noindent In 1955, Amitsur \cite{Amitsur55} classified the finite subgroups of the multiplicative group of a division algebra whose central field $k$ has characteristic zero. This completed work of Herstein \cite{Herstein53} who treated the case the central field had prime characteristic. The important class of those finite subgroups which do arise in division algebras in regard to this article are the \emph{$D$--groups} defined by
\[ G_{m,r} = \innp{X,Y~:~X^m=1,~Y^n = X^t,~YXY^{-1} = X^r}, \]
where $(m,r)=1$, $r<m$, $s=(r-1,m)$, $m=st$, and $n$ is the multiplicative order of $r$ in $\Z/m\Z$. The importance of $D$--groups is seen in the following result---this is a special case of Corollary 7 in \cite{Amitsur55}.\smallskip\smallskip

\begin{thm}[Amitsur]\label{T:Amitsur} If $A$ is a cyclic $E$--division algebra of odd degree $p$ with $\ch(E)=0$ and $G$ is a finite subgroup of $A^\times$, then $G=G_{m,r}$ and $n\mid p$.
\end{thm}

\noindent It is worth noting Amitsur's paper provides even more information than what we have stated here. It also holds when $A$ is not necessarily cyclic, a case which can never happen when $E$ is a number field.

%---------------------------------------------------------------------------
\subsection{Proof of \refT{PCST}}

\noindent We now commence with the proof of \refT{PCST}.

\begin{pf}[Proof of \refT{PCST}] Let $G$ be a finite subgroup of an arithmetic lattice $\Ga$ of second type. By construction, $\Ga$ is contained in $A^\times$ for a cyclic central division algebra $A$ of degree $p$ over a CM field $E/F$. In particular, $G$ is a subgroup of $A^\times$ and thus \refT{Amitsur} implies $G \cong G_{m,r}$ with $n=1$ or $p$. In the former case, $G_{m,r}$ is cyclic of order $m$. Otherwise, we have the short exact sequence
\[ 1 \lra \innp{X} \lra G_{m,r} \lra \innp{\ol{Y}} \lra 1 \]
where $\innp{X}$ is the normal cyclic subgroup of $G_{m,r}$ of order $m$ and $\innp{\ol{Y}}$ is the cyclic subgroup of order $p$ generated by the image of $Y$ under the canonical projection. As $\Ga$ is a lattice in $\Uni(p-1,1)$, we have a faithful representation
\[ \rho\co \Ga \lra \Uni(p-1,1). \]
Restricting $\rho$ to $G$ produces the faithful representation
\[ \rho\co G \lra \Uni(p-1,1). \]
As a finite subgroup of $\Uni(p-1,1)$, $\rho(G)$ is contained in a maximal compact subgroup of $\Uni(p-1,1)$. However, any maximal compact subgroup of $\Uni(p-1,1)$ is conjugate in $\Uni(p-1,1)$ to $\Uni(p-1) \times \Uni(1)$ (\cite[Ch. 1]{Knapp96}). Therefore, the representation $\rho$ is reducible and decomposes into a nontrivial direct sum of representations
\[ \rho = \bop_{j=1}^\ell \rho_j, \quad p = \sum_{j=1}^\ell d_j, \quad d_j = \dim \rho_j. \]
In particular, the reducibility of $\rho$ yields the inequality $d_j < p$ for all $j=1,\dots,\ell$. The irreducible representations of the groups $G_{m,r}$ are given by induction on the irreducible representations of $\innp{X}$ (\cite[p. 61]{Serre77}) and have degrees which divides $[G_{m,r}:\innp{X}]=p$. As each $d_j<p$, and $d_j \mid p$, the primality of $p$ implies each $d_j$ must be one. Hence, $\rho$ is a direct sum of 1--dimensional representations. However, $G_{m,r}$ is nonabelian for $n=p$ and so cannot possess a faithful representation whose summands are all 1--dimensional. Hence $G$ must be cyclic.
\end{pf}

\noindent The proof of \refT{PCST} works more generally for arithmetic lattices of second type in the groups $\Uni(p_1,p_2)$ for $p_1+p_2=p$ ($p_1,p_2>0$). The lone point of deviation between this setting and $\Uni(p-1,1)$ is the $\Uni(p_1,p_2)$--conjugacy of any finite subgroup into the maximal subgroup $\Uni(p_1)\times \Uni(p_2)$ opposed to $\Uni(p-1) \times \Uni(1)$. However, this conjugacy's solitary use was in the deduction of the reducibility of the representation $\rho$, a fact which still visibly persists.\smallskip\smallskip

\begin{cor} If $\Ga$ is an arithmetic lattice of second type in $\Uni(p_1,p_2)$ with $p_1+p_2=p$ for an odd prime $p$ and $p_1,p_2>0$, and $G$ is a finite subgroup of $\Ga$, then $G$ is cyclic.
\end{cor}

%---------------------------------------------------------------------------
%---------------------------------------------------------------------------
\section{Concluding remarks}\label{S:Final}

\noindent More generally one can ask which finite groups arise in lattices of other Lie groups. For arithmetic lattices in the real classical groups, we will address this question in a future paper. For many classes of lattices, the results and proofs are the same as the ones given here. For instance, the arithmetic lattices arising from bilinear and hermitian forms in $\SO(p,q)$ and $\SU(p,q)$ have existence theorems identical to \refT{AllFiniteAllClasses}. The arithmetic lattices in $\Sp(p,q)$ also possess this existence property as every arithmetic lattice in $\SO(p,q)$ or $\SU(p,q)$ arising from a form injects into an arithmetic lattice of $\Sp(p,q)$. The classes of arithmetic lattices of most interest in regards to this problem are:
\begin{itemize}
\item Arithmetic lattices in $\SU(p,q)$ of second or mixed type.
\item Arithmetic lattices in $\SO(p,q)$ arising from quaternion algebras.
\item Arithmetic lattices in $\SL(n;\R)$ and $\SL(n;\C)$.
\end{itemize}
The lattices in $\SU(p,q)$ of mixed type can be regarded as intermediate or transition lattices between first and second type. As such, the finite subgroups they possess fall in between those possessed by first type lattices and those possessed by second type lattices. 

%---------------------------------------------------------------------------
%---------------------------------------------------------------------------

\def\cprime{$'$} \def\lfhook#1{\setbox0=\hbox{#1}{\ooalign{\hidewidth
  \lower1.5ex\hbox{'}\hidewidth\crcr\unhbox0}}} \def\cprime{$'$}
  \def\cprime{$'$}
\providecommand{\bysame}{\leavevmode\hbox to3em{\hrulefill}\thinspace}
\providecommand{\MR}{\relax\ifhmode\unskip\space\fi MR }
% \MRhref is called by the amsart/book/proc definition of \MR.
\providecommand{\MRhref}[2]{%
  \href{http://www.ams.org/mathscinet-getitem?mr=#1}{#2}
}
\providecommand{\href}[2]{#2}

%---------------------------------------------------------------------------
%---------------------------------------------------------------------------

\noindent
Department of Mathematics, The University of Texas at Austin\\
email: {\tt dmcreyn@math.utexas.edu}

%---------------------------------------------------------------------------
%---------------------------------------------------------------------------

\end{document}